\documentclass[12pt]{article}
\usepackage[utf8]{inputenc}
\usepackage{amsmath, amssymb, url, float, graphicx, float}
\usepackage{fullpage}
\usepackage[small,bf]{caption}
\usepackage{cite}
\usepackage{xcolor}
\usepackage{subcaption}

\renewcommand\phi\varphi

\let\phi\varphi

\newcommand{\ones}{\mathbf 1}
\newcommand{\reals}{{\mbox{\bf R}}}

\newcommand{\eg}{{\it e.g.}}
\newcommand{\ie}{{\it i.e.}}

\newcommand{\BEAS}{\begin{eqnarray*}}
\newcommand{\EEAS}{\end{eqnarray*}}
\newcommand{\BEA}{\begin{eqnarray}}
\newcommand{\EEA}{\end{eqnarray}}
\newcommand{\BEQ}{\begin{equation}}
\newcommand{\EEQ}{\end{equation}}
\newcommand{\BIT}{\begin{itemize}}
\newcommand{\EIT}{\end{itemize}}

\newif\iftodos

\title{A Note on the Welfare Gap in Fair Ordering}
\author{
    Theo Diamandis\\
    {\small \texttt{tdiamand@mit.edu}} 
    \and 
    Guillermo Angeris\\
    {\small \texttt{gangeris@baincapital.com}}
}
\date{March 2023}

\begin{document} 
\maketitle 

\begin{abstract}
    Public blockchains group submitted transactions into batches, called blocks.
    A natural question is how to determine 
    which transactions are included in these batches.
    In this note, we show a gap between the welfare of so-called
    `fair' ordering, namely first-in-first-out (an ideal that a number of
    blockchain protocols strive to achieve), where the first transactions to
    arrive are the ones put into the block, and the welfare of `optimal' inclusion that is,
    at least approximately, welfare-maximizing, such as choosing which
    transactions are included in a block via an auction. We show this gap is
    positive under a simple model with mild assumptions where we assume
    transactions are, roughly speaking, uniformly drawn from a reasonable
    distribution.
    Our results formalize a performance metric for blockchain inclusion rules
    and consequently provide a framework to help design and compare these rules.
    The results can be directly extended to ordering mechanisms as well.
\end{abstract}

\section*{Introduction}
In this note, we examine ways in which first-in-first-out blockchains
distort their block space markets, and we quantify the resulting decrease in social
welfare. While we focus on the blockchain setting,
these results are directly applicable to any batched system that must allocate
its finite resources among different users.

Public blockchains generally implement a \emph{fee mechanism} to allocate finite
computational resources across competing transactions. These transactions are
broadcasted by users through the peer-to-peer network and collected in the
mempool, which contains submitted transactions that have not yet been included
on chain. A validator then selects which transactions from the mempool are 
included in the next block, possibly subject to an \emph{inclusion rule}. (We
note that this view is overly simplistic, but it suffices for the purposes of 
this note.) 

Most protocols implement either a fixed fee per transaction, or a dynamic fee
that fluctuates with the demand for block space (\eg, the base fee on
Ethereum). The majority of commonly-used protocols implement an unconstrained
inclusion rule: validators can include and order transactions in a block
however they wish. However, many authors have proposed a first-in-first-out
(FIFO) inclusion rule (also called first-come-first-served or `fair' ordering), 
where transactions are included in the order in which
they arrive in the mempool, until the block is full. These
methods (for example, those proposed in~\cite{aequitas, themis,pompe,
quickorderfairness, wendy}) add extra rules to consensus to ensure, roughly
speaking, that if a majority of validators receive transaction $A$ before
transaction $B$, then transaction $A$ is guaranteed to be included in the
block, if transaction $B$ is included. (We assume this is achievable and
ignore implementation difficulties of FIFO inclusion in this note.) 

\paragraph{Externalities.} The simplest case to consider is where the
transaction fee is fixed and the transaction inclusion rule is (effectively)
FIFO. The fixed---usually very low---transaction fee, coupled with FIFO
inclusion, encourages users competing for specific transactions (\eg, arbitrage
trades or liquidations) to spam the network in hopes that their transaction is
seen first and therefore included in the block. This spam may result in block
space filled with reverting, `garbage' transactions with low utility. A recent
analysis of Solana transactions found that over 50\% of transactions were
failed arbitrage trades~\cite{jito2023}. Clearly, these users are not paying
for the externality they are causing to the network, negatively impacting other
users. A natural solution to spam is to charge a dynamic fee for
block space: as the demand increases, the fee increases (see
EIP-1559~\cite{eip1559}). In this note we will see that, even in this case,
FIFO inclusion still can dramatically reduce the social welfare generated by
the network. As a simple example, liquidations (which are high-utility
transactions) would not be prioritized, again encouraging users to spam the
network to ensure their transaction is included rather than simply paying a
higher fee. This type of challenge played out in the Arbitrum ecosystem:
although there is a base fee, FIFO inclusion still encouraged spamming the
network for transaction prioritization~\cite{snoopy2023}. We consider this case
for the remainder of this note.

\section{Building blocks}
In this section, we briefly introduce the block building problem, which is a
special case of the 0-1 knapsack problem~\cite{kellerer2004knapsack}. We then
show a simple lower bound to the optimal value of this problem via a greedy
heuristic, which can be computed efficiently. We will use this lower bound in a
later section to bound the welfare gap between the `FIFO' inclusion and the best
possible transactions to include (to maximize social welfare, \ie, net utility).

\paragraph{The block building problem.} Following~\cite{diamandis2022dynamic},
we consider a mempool with transactions $j=1, \dots, n$. Each transaction $j$
has utility $\tilde q_j \in \reals_+$ and consumes some amount of a single 
resource, called \emph{gas}, denoted by $a_j \in \reals_+$. Gas has a per-unit 
cost $g \in \reals_+$ in the same num\'eraire as the utility. We denote the gas 
limit per block by $b \in \reals_+$. 
We define the \emph{block building problem} as
\begin{equation}
    \label{eq:block-building}
    \begin{aligned}
        & \text{maximize} && q^Tx\\
        & \text{subject to} && a^Tx \le b  \\
        &&& x \in \{ 0,1 \}^n,
    \end{aligned}
\end{equation}
where $q = \tilde q - ga$ is the vector of net transaction utilities, 
after removing gas fees. Here,
the variable $x \in \reals^n$ indicates which transactions are included in the
next block: $x_j = 1$ if transaction $j$ is included in the block and $x_j
= 0$ otherwise. In contrast to~\cite{diamandis2022dynamic}, which considers arbitrarily
complicated constraints, we only consider the constraint that the gas 
total of all included transactions must be
lower than the upper limit. We assume that $q$ is nonnegative (\ie, users do
not submit negative utility transactions) and 
denote the optimal value of problem~\eqref{eq:block-building} by $p^\star$.

\paragraph{Discussion.} Problem~\eqref{eq:block-building} is an instance of the
weighted 0-1 knapsack problem~\cite{kellerer2004knapsack}. There
are a number of important issues about how, exactly, to instantiate this
problem and solve it in practice. The first: how does one elicit the
utilities $q$ from users? This vector $q$ can come from certain mechanisms
such as second-price auctions~\cite{roughgarden2016twenty}, though
implementing these on chain remains an open research problem
(see~\cite{chitra2023credible,MEV_Suave, chung2021foundations} and references
therein). We do not deal with this question here and assume that there
is some mechanism for receiving (or approximating) these utilities. 
Additionally, this problem is NP-hard to solve in the worst case~\cite[App.\ A]{kellerer2004knapsack}. 
However, in practice, problems of this form are commonly solved to
optimality using software such as Gurobi~\cite{gurobi}, and good heuristic 
solutions can be efficiently computed. 

\paragraph{FIFO inclusion.} In FIFO inclusion, on the other hand, transactions
are included `as they arrive', until the block reaches capacity. We assume that
the transactions $j = 1, \dots, n$ are ordered by arrival time, from earliest
to latest. (Whether FIFO inclusion is exactly possible in a decentralized
setting has been a subject of great debate; see, for example,~\cite[Theorem 1.2]{kelkar2020order}
and~\cite{vafadar2023condorcet}.
However, we assume the `ideal' scenario that many protocols are trying to
achieve.) In this scenario, the vector $x$ must have the form $x = (\ones_k,
0)$ for some $k$. We can define the utility from the FIFO transaction
explicitly as the maximum utility over over all vectors of this form:
\[
    p^\mathrm{FIFO} = \max\{q^T(\ones_k, 0) \mid a^T(\ones_k, 0) \le b,\; k = 0, 1, \dots, n \}.
\]
We aim to answer the question `what is difference between the total utility of an
optimal block packing (\ie, a solution to problem~\eqref{eq:block-building})
and that of a FIFO packing?'. In other words, we want to find a lower bound on the
gap
\[
    \Gamma = p^\star - p^\mathrm{FIFO}.
\] 
Of course $\Gamma \ge 0$, since any FIFO solution is feasible for
the original problem~\eqref{eq:block-building}. Though, we will show that
the quantity $\Gamma$ is positive (and potentially large) `in expectation'
under weak assumptions about the distribution of the transactions.

\subsection{A heuristic for block building}\label{sec:heuristic}
Since problem~\eqref{eq:block-building} is a weighted knapsack problem, a number
of simple heuristics provide approximate solutions with relatively
tight guarantees. We present a very basic overview of one important heuristic
and its corresponding proof of tightness here. We will then use this heuristic,
and corresponding bound, to approximate $p^\star$ and show that $\Gamma$ is large
in a number of common scenarios.

\paragraph{Greedy heuristic.} The simplest (and most common) heuristic is to
relax the integrality constraint ($x \in \{0, 1\}^n$) in
problem~\eqref{eq:block-building} to an interval constraint, to
get the linear program:
\begin{equation}\label{eq:relaxation}
    \begin{aligned}
        & \text{maximize} && q^Tx\\
        & \text{subject to} && a^Tx \le b\\
        &&& 0 \le x \le \ones,
    \end{aligned}
\end{equation}
with variable $x \in \reals^n$ and the same problem data as~\eqref{eq:block-building}.
We write the optimal value of this relaxed problem as $r^\star$. Note that, since
every $x$ that is feasible for~\eqref{eq:block-building} is feasible for its
relaxation~\eqref{eq:relaxation}, we have that 
\[
    r^\star \ge p^\star \ge p^\mathrm{FIFO}.
\]
Problem~\eqref{eq:relaxation} is easy to (computationally) solve in practice. In fact, it is
possible to write a closed-form solution $x^\star$ to~\eqref{eq:relaxation}. 
To see this, start with the equivalent problem,
\[
    \begin{aligned}
        & \text{maximize} && \bar q^Ty\\
        & \text{subject to} && \ones^Ty \le 1\\
        &&& 0 \le y_i \le a_i/b \quad i=1, \dots, n.
    \end{aligned}
\]
Here, the problem data are the \emph{efficiencies} $\bar q_i = bq_i/a_i$ for
each transaction $i=1, \dots, n$, while the rest of the definitions are
identical to the original relaxation~\eqref{eq:relaxation}. In this problem, we
have simply done a variable substitution $y_i = a_ix_i/b$ in the relaxed
problem~\eqref{eq:relaxation}. We can think of the $\bar q_i$ as a measure of
the `welfare-per-unit-resource' for transaction $i$.

A solution to this problem is very simple to construct: sort the $\bar
q_i$ in nonincreasing order, with indices $\tau_1, \dots, \tau_n$, then
set $y_{\tau_i}^\star = a_{\tau_i} / b$, in order, for each $i$ until the
constraint $\ones^T y \le 1$ is met. If there is no
index at which the constraint is met, an optimal solution is to set
$y_{\tau_i}^\star = a_{\tau_i} / b$ for $i=1, \dots, n$. (It is \emph{the}
optimal solution if $\bar q > 0$.) If there is an index at which the constraint is
met or surpassed, say index $\tau_k$, then we choose $y_{\tau_k}^\star$ to be
the largest possible value less than or equal to $a_{\tau_k} / b$ such that the 
constraint is met at
equality. (There may be many solutions if there are many entries of $\bar q$
with value $\bar q_{\tau_k}$, in which case any entry suffices.) We can
recover a solution $x^\star$ for problem~\eqref{eq:relaxation} from this
optimal $y^\star$ by using the substitution above, $x_i^\star =
by_i^\star/a_i$ for $i=1, \dots, n$. Note that this solution, $x^\star$, is
fractional (\ie, has $0 < x_i^\star < 1$ for some $i$) in at most one entry.

\paragraph{Discussion.} The solution to~\eqref{eq:relaxation} above suggests 
a reasonable heuristic to solve~\eqref{eq:block-building}: since at most one 
entry of an optimal solution $x^\star$ to~\eqref{eq:relaxation} is fractional, we can simply round
the (at most one) fractional entry down to zero, to produce $x^0 \in \{0,
1\}^n$. Note that $x^0$ is feasible for the original
problem~\eqref{eq:block-building}, and we denote its objective value as $p^0 =
q^Tx^0$. By definition, we have that $p^0 \le p^\star$, since $x^0$ is
feasible. Perhaps surprisingly, we also have that $p^0 \ge Cp^\star$ for some constant
$C \le 1$. In other words, this heuristic gives a feasible point with objective
value $p^0$, which is `close' to the true optimal value $p^\star$.

\paragraph{Bound.} If the transactions all consume gas at most $a_i \le b/m$ for some
integer $m > 1$, then we have that
\[
    \frac{m}{m-1}p^0 \ge p^\star \ge p^0.
\]
In other words, the heuristic solution, with optimal value $p^0$, is very close
to the optimal value of exactly solving problem~\eqref{eq:block-building},
which is NP-hard, whenever $m$ is somewhat large. We expect this to often be
the case in practice: high-value transactions like liquidations and arbitrages
generally don't consume much gas, relative to the block limit.

\paragraph{Proof.}
To see this bound, first note that the solution to the relaxation~\eqref{eq:relaxation}
satisfies $a^Tx^\star = b$. By definition, there is at most one
nonintegral entry, with gas $a_i \le b/m$ so the rounded solution must satisfy
$a^Tx^0 \ge (m-1)b/m$. Let $i$ be the entry with largest $q_i$ such that $x_i^0 = 0$
(\ie, let $i$ be the highest-utility transaction not included in the heuristic solution),
then
\[
    \frac{q_i}{a_i} \le \frac{q^Tx^0}{a^Tx^0} = \frac{p^0}{a^Tx^0} \le \frac{m}{(m-1)b}p^0.
\]
The first inequality follows from the fact that, for
positive $t, u, v, w$ we have
\[
    \min\left\{\frac{t}{v}, \frac{u}{w}\right\} 
    \le \left(\frac{v}{v+w}\right)\frac{t}{v} + \left(\frac{w}{v+w}\right)\frac{u}{w}
    \le \frac{t + u}{v + w},
\]
and $q_i/a_i$ is no larger than the entries $j$ with $x_j^0 = 1$, by construction of $x^0$.
The equality follows by definition, and the last inequality comes from the
previous discussion.
Finally, using the fact that $a_i \le b/m$ again and the above, we have that
\[
    \frac{m}{b}q_i \le \frac{q_i}{a_i} \le \frac{m}{(m-1)b}p^0,
\]
or, simplifying,
\begin{equation}\label{eq:bound-q}
    q_i \le \frac{p^0}{m-1}.
\end{equation}
Since there is at most one transaction partially included by the relaxation,
and this transaction has utility no larger than $q_i$ (by definition of $i$),
then
\[
    p^0 + q_i \ge r^\star \ge p^\star \ge p^0,
\]
where, from before, $r^\star$ is the optimal value of the
relaxation~\eqref{eq:relaxation}. Combining this with~\eqref{eq:bound-q}, we then have
that
\[
    \frac{m}{m-1}p^0 \ge p^\star \ge p^0.
\]
For more on similar approximations, see~\cite[\S6]{kellerer2004knapsack}.

\section{What's the gap?}
Given the discussion above, it makes sense to consider a model where transactions
are drawn from some distribution, which we will characterize in terms of the
efficiencies. We will also assume that there are minimum and maximum gas limits 
for any transaction, denoted by $B^-$ and $B^+$ respectively, so that 
$B^- \le a_i \le B^+$ for $i = 1, \dots, n$. Clearly, the minimum and maximum
transaction sizes imply a maximum and minimum number of transactions we can
include in each block, given by $b/B^-$ and $b/B^+$, respectively. From this assumption,
we can construct a lower bound of the optimal block utility and an upper bound 
of the FIFO-inclusion block utility. In this section, we will construct these
bounds and discuss when there is a strictly positive gap between the utilities.

\paragraph{Lower bound.}
First, we will lower bound the utility of a block that was packed using the greedy 
heuristic (which we know is close to optimal).
Denote the number of transactions included in the block using the greedy
heuristic by $\bar k$. We define the average utility of these transactions by
\[
    q^+ = \frac{1}{\bar k}\sum_{i=1}^{\bar k} q_{\tau_i}.
\]
Since $\bar k \ge b/B^+$, 
a lower bound for the utility of the greedily packed block is
\[
    L = \bar k q^+  \ge \frac{b}{B^+}q^+.
\]

\paragraph{Upper bound.}
Now we will upper bound the expected utility of a FIFO-inclusion block, assuming the
arrival time of the transactions is random, \ie, that the $n$
transactions are uniformly randomly permuted. We know at
most $b/B^-$ transactions can be included in a block, which means that the expected
utility of FIFO, which we will call $U$, is no larger than
\[
    U \le \frac{b}{B^-}\frac{1}{n}\sum_{i=1}^n q_i.
\]
If we define $q^-$ as the average utility for the transactions not included by
the greedy heuristic,
\[
    q^- = \frac{1}{n-\bar k}\sum_{i=\bar k+1}^{n} q_{\tau_i},
\]
then we can write
\[
    \frac{1}{n}\sum_{i=1}^n q_i = \frac{\bar k}{n}q^+ + \left(1-\frac{\bar k}{n}\right)q^-.
\]
(We can view $q^-$ as, roughly speaking, the average utility of the `tail' of
transactions, as the efficiencies get small.) This means the average FIFO
block utility is bounded from above by
\[
    U \le \frac{b}{B^-}\left(\frac{\bar k}{n}q^+ + \left(1-\frac{\bar k}{n}\right)q^-\right).
\]

\paragraph{What's the gap?}
We now characterize the gap between these two bounds, which
gives us a lower bound on the utility gap between the FIFO inclusion block and 
the optimal block, since
\[
    \Gamma \ge L - U \ge \frac{b}{B^+}q^+ - \frac{b}{B^-}\left((q^+ - q^-)\frac{\bar k}{n} + q^- \right).
\]
We give a basic condition for when the right-hand-side of this
inequality is positive, which would imply that the gap $\Gamma > 0$. Rearranging,
it is easy to see that the right hand side is positive whenever
\begin{equation}\label{eq:main-bound}
    q^+\left( 1 - \frac{\bar k \eta}{n}\right) >   
    \eta q^-\left(1 - \frac{\bar k}{n}\right),
\end{equation}
where $\eta = B^+ / B^-$ is the ratio between the largest and smallest possible
transaction. If
    $q^+ > \eta q^-$,
then, as the number of outstanding transactions $n$ becomes large relative to
the number of greedily-chosen transactions, $\bar k$, while the average
utilities stay roughly constant, we get that $\Gamma > 0$.
We may also wish to consider the ratio of optimal to FIFO block utility. From
the preceding discussion, we have that
\[
    \frac{p^\star}{p^\mathrm{FIFO}} \ge \frac{L}{U} \ge
    \frac{q^+/\eta}{(q^+ - q^-)\frac{\bar k}{n} + q^-},
\]
which is strictly greater than one under the positive gap condition~\eqref{eq:main-bound}.

\paragraph{Discussion.}
We note that the bound derived in~\eqref{eq:main-bound} is actually loose 
in practice due to the worst-case
assumptions made in the lower bound, where we only use the fact that $\bar k
\ge b/B^+$ and the upper bound, where we use the fact that $b/B^-$ bounds the
maximum possible number of transactions that can be included by FIFO. Both of
these bounds are quite loose if there are only a few transactions close to
the gas bounds, relative to the rest of the distribution. In this note,
we only seek `qualitatively reasonable' behavior from the bounds in order to
gain intuition; we don't seek tight constants. An interesting avenue for future
research would be to either tighten the bounds given here or to generalize them to
the multidimensional setting considered in~\cite{diamandis2022dynamic}.

\subsection{Special cases and experiments} 
From our bound~\eqref{eq:main-bound}, we can
deduce some conditions under which we are guaranteed to have a large gap
between the utilities of the optimal and FIFO-inclusion blocks. We also show 
basic numerical experiments which suggest that, in practice, the gap may be 
far larger than the one suggested by our bounds.

\paragraph{Sharp distribution.} If the distribution over the
efficiencies is sufficiently sharp (\eg, if there are a small number of transactions
with very large utility and equal or lower gas relative to the others) then
$q^+ \gg q^-$ for reasonable block sizes. This situation is common in practice
when many similar transactions are submitted but only one can be executed
profitably, as is the case in many MEV opportunities, including liquidations and
DEX arbitrage. 
Since only one of these transactions can be executed
profitably, the others will revert and, therefore, have non-positive utility. 
Note that this situation corresponds to these efficiencies being drawn from a
heavy-tailed probability distribution.

\paragraph{Equal size transactions.} If the transactions have roughly the 
same size, \ie, $\eta \approx 1$, then our bound~\eqref{eq:main-bound} implies 
that any distribution of utilities that is
not flat will create a gap between the optimal and FIFO-inclusion blocks.
Intuitively, this is straightforward to see: there is always a chance that
high-utility transactions will not be included in the block, and, since all
transactions consume roughly the same gas, this is a strict loss in total
utility.

\begin{figure}[h!]
    \captionsetup[sub]{font=scriptsize}
    \centering
    \begin{subfigure}[t]{0.32\linewidth}
        \centering
        \includegraphics[width=\columnwidth]{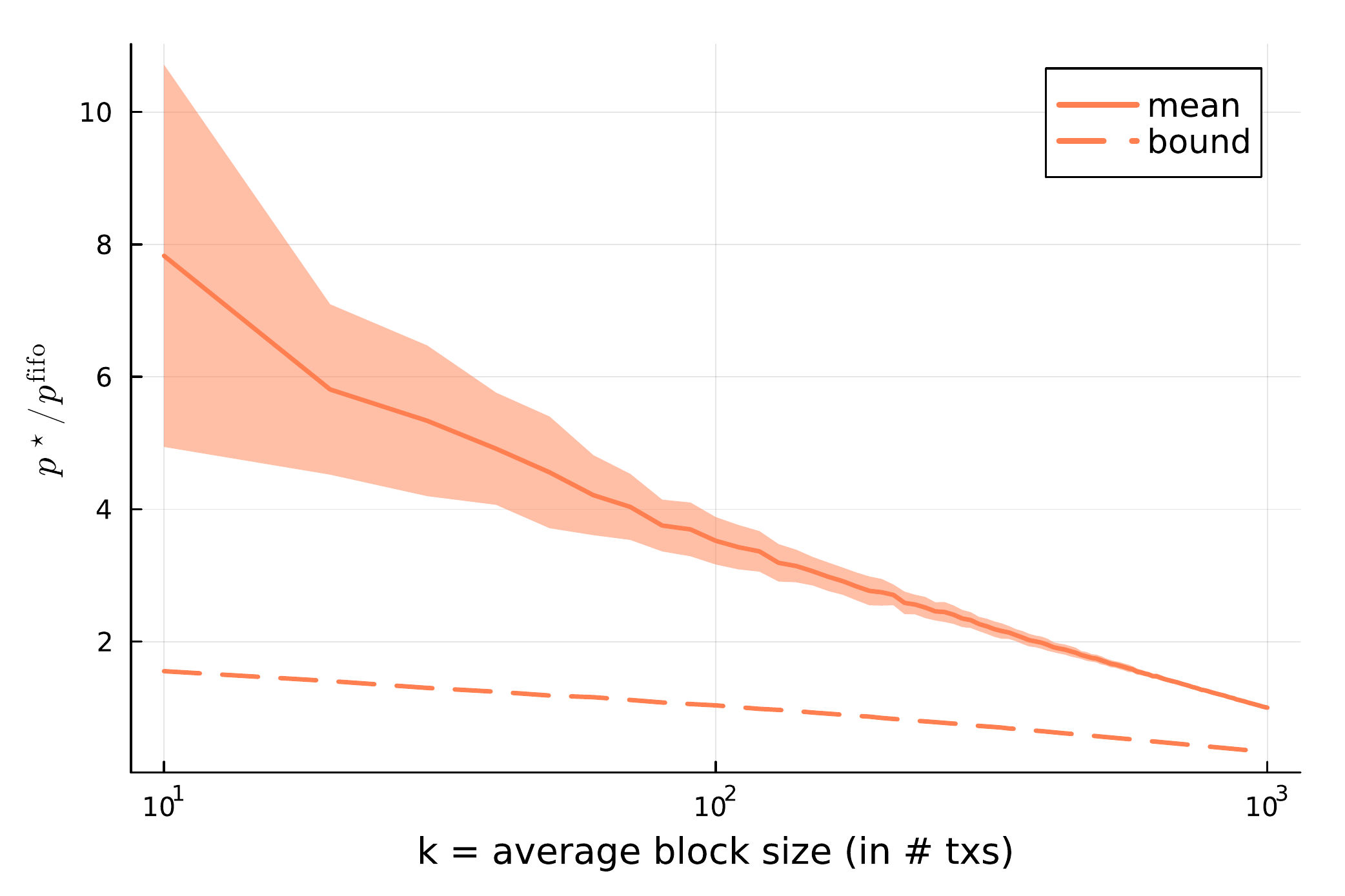}
        \caption{$q_i \sim \mathrm{Exponential}(2.5)$}
    \end{subfigure}
    \hfill
    \begin{subfigure}[t]{0.32\linewidth}
        \centering
        \includegraphics[width=\columnwidth]{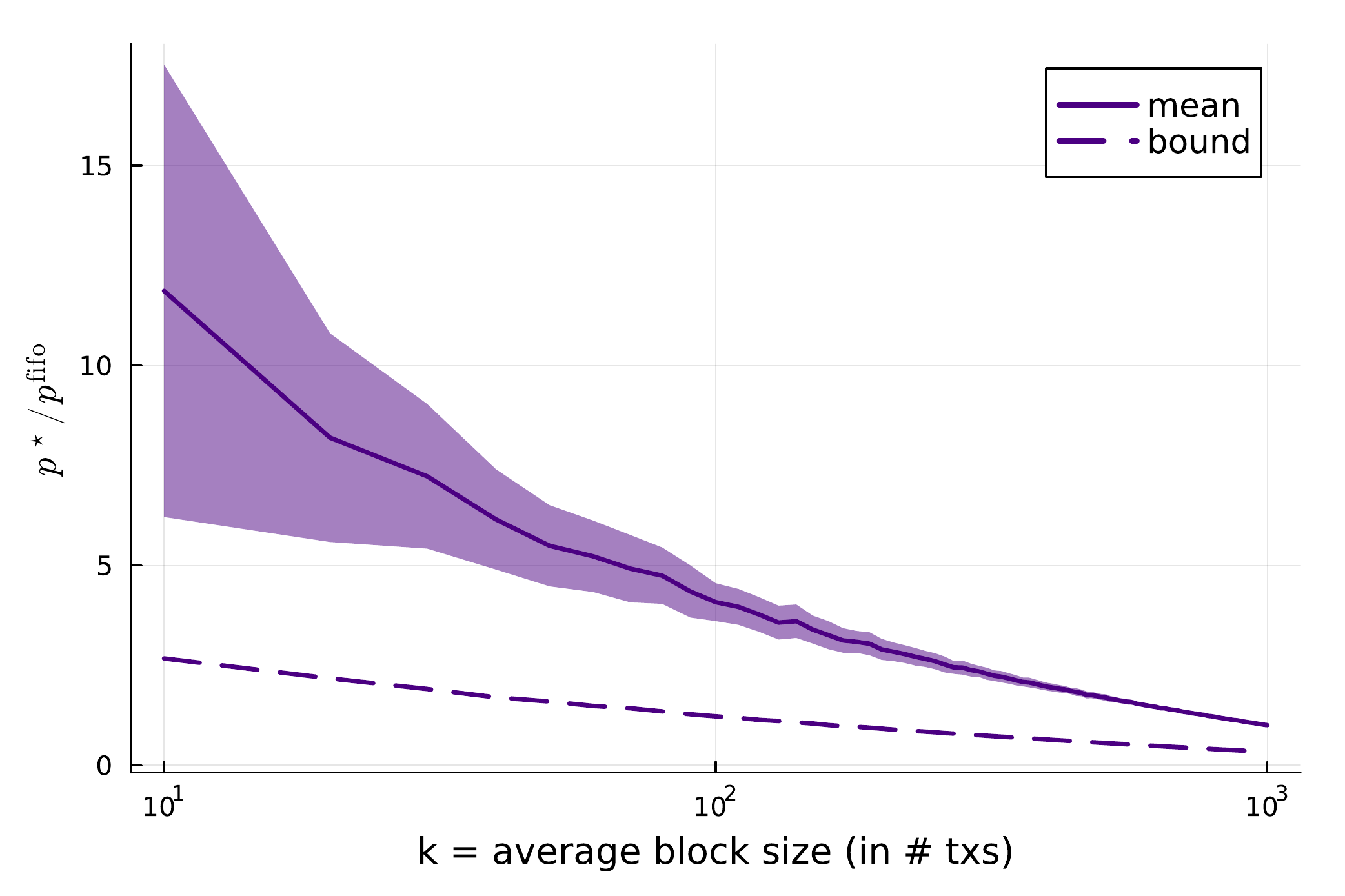}
        \caption{$q_i \sim \mathrm{LogNormal}(1, 1)$}
    \end{subfigure}
    \hfill
    \begin{subfigure}[t]{0.32\linewidth}
        \centering
        \includegraphics[width=\columnwidth]{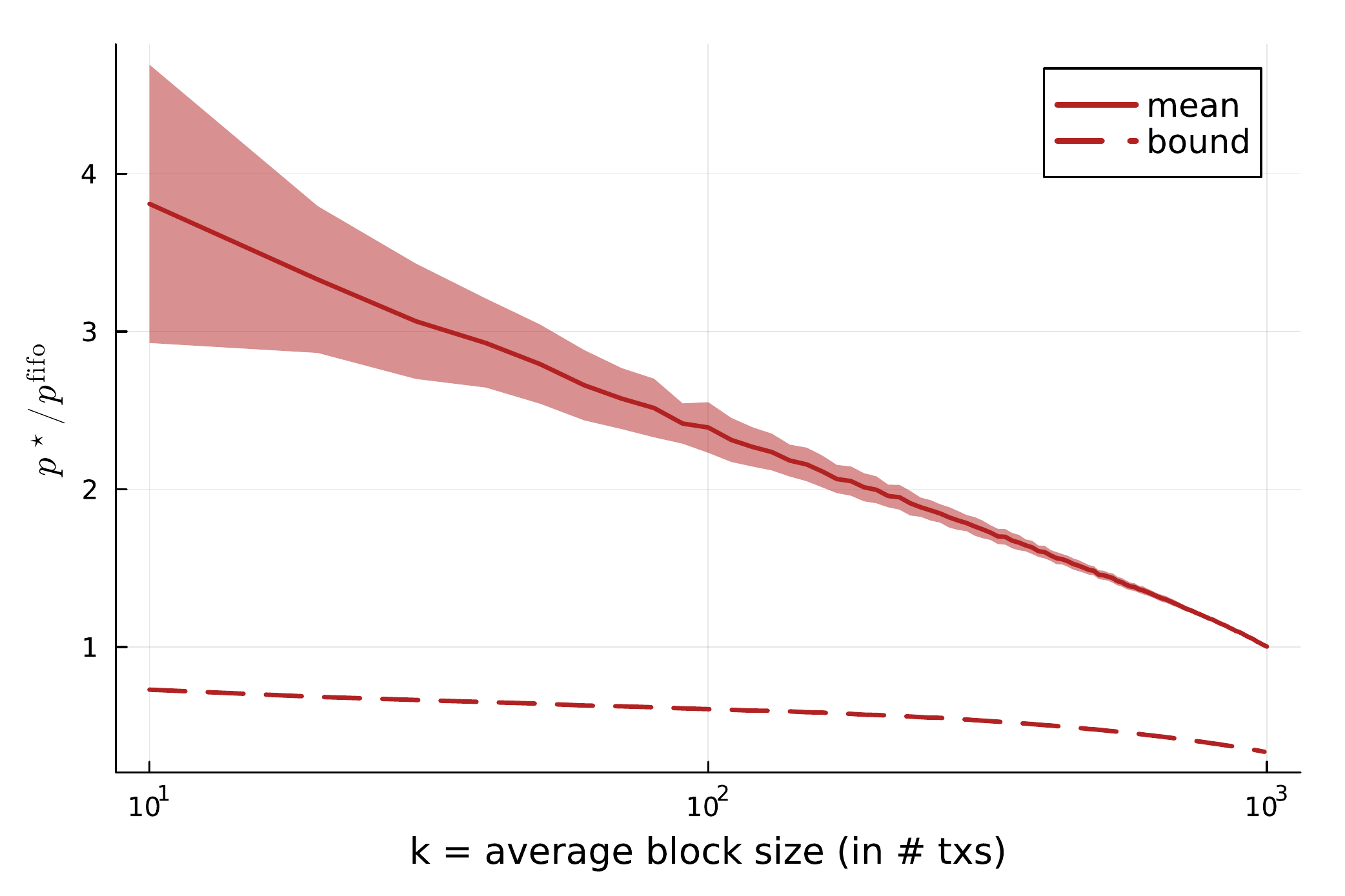}
        \caption{$q_i \sim \mathrm{Rayleigh}(1)$}
    \end{subfigure}
    \caption{The gap between FIFO and optimal block packing for flat 
    (light-tailed) utility distributions is strictly positive but decreases as 
    block size increases.
        }
    \label{fig:gap}
\end{figure}

\begin{figure}[h!]
    \captionsetup[sub]{font=scriptsize}
    \centering
    \begin{subfigure}[t]{0.45\linewidth}
        \centering
        \includegraphics[width=\columnwidth]{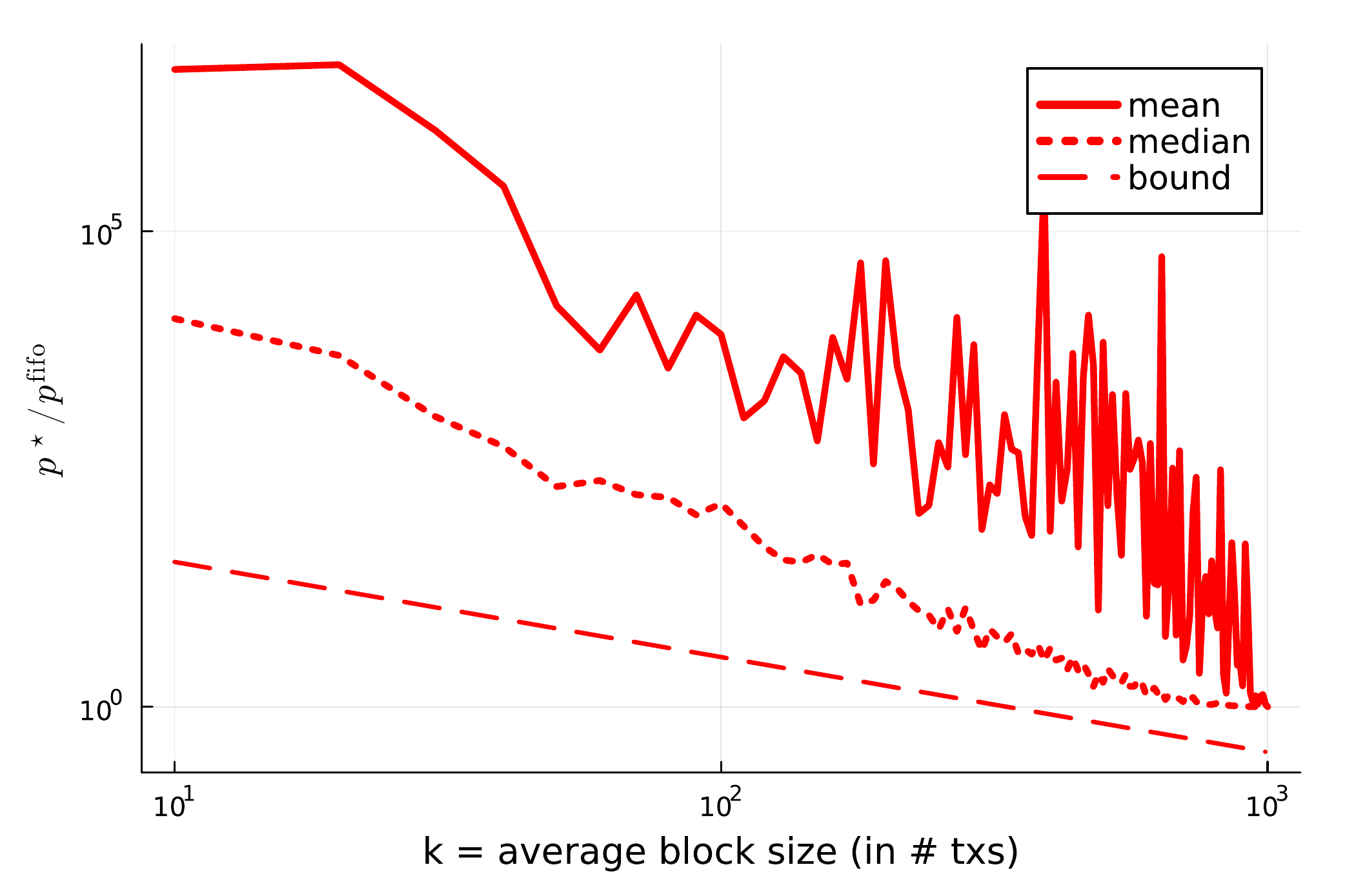}
        \caption{$q_i \sim \mathrm{Levy}(0, 1)$}
    \end{subfigure}
    \hfill
    \begin{subfigure}[t]{0.45\linewidth}
        \centering
        \includegraphics[width=\columnwidth]{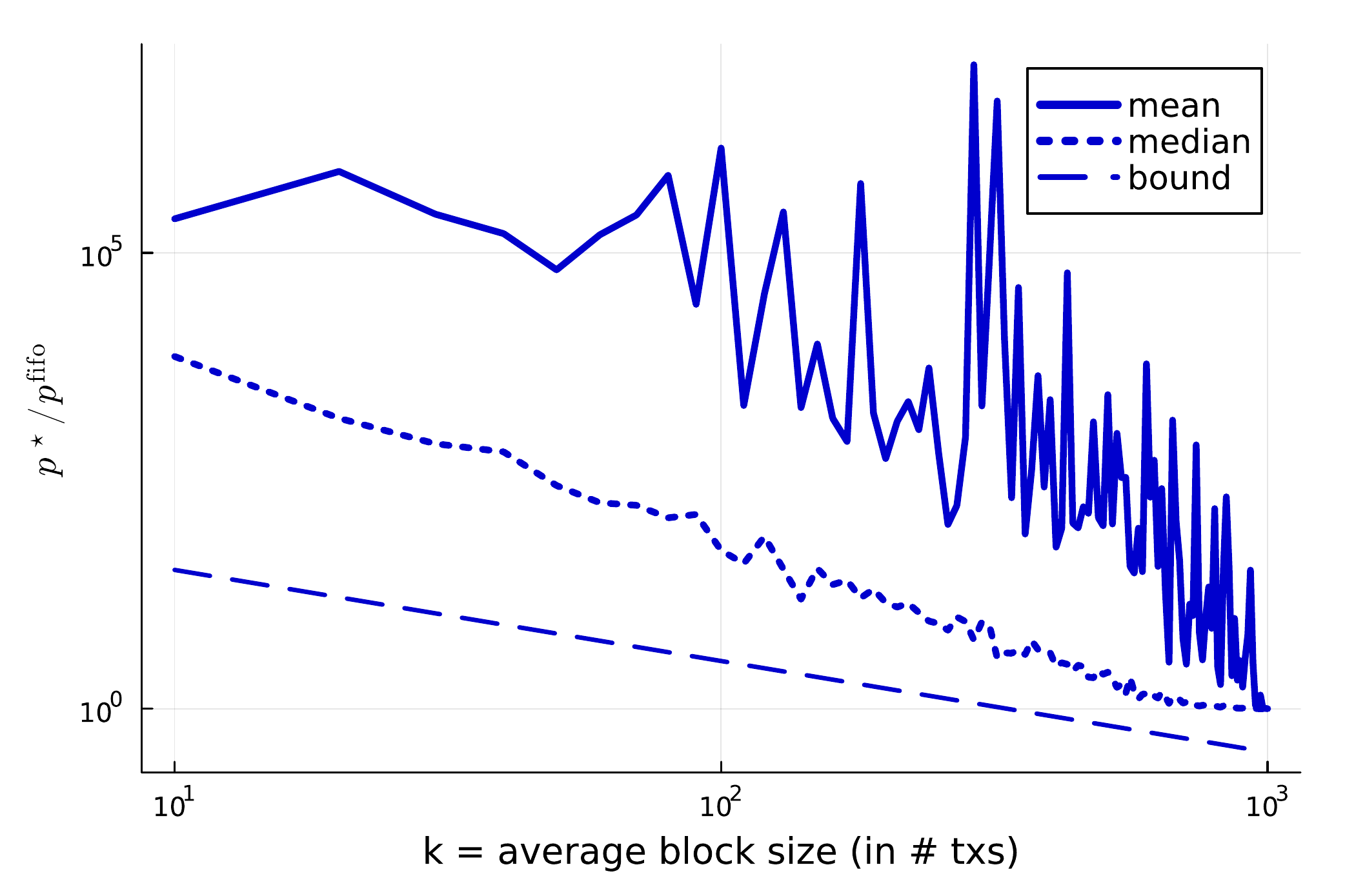}
        \caption{$q_i \sim \mathrm{Pareto}(0.5)$}
    \end{subfigure}
    \caption{The gap between FIFO and optimal block packing for sharp 
    (heavy-tailed) distributions, which are common in practice, is substantially
    larger.
    }
    \label{fig:gaps-heavy}
\end{figure}

\paragraph{Simple experiments.}
We plot our utility gap ratio bound and the realized 
gap for several transaction utility distributions
(we provide definitions in appendix~\ref{app:dists}).
We draw transaction sizes uniformly at random from the interval $[1, 3]$, 
generate $1000$ transactions in the mempool, and vary the block size from $20$
to $2000$ gas (for an average of $10$ to all $1000$ transactions per block).
We run $100$ trials for each block size.
All experiments use the Julia programming language~\cite{bezanson2017julia}
and the \texttt{Distributions.jl} package~\cite{distributions.jl}.
Code is available at
\begin{center}
    \texttt{https://github.com/bcc-research/fifo-note}
\end{center}
Figure~\ref{fig:gap} shows the empirical gaps for flat distributions (\ie, those where we do not
expect transactions with very large utility, relative to the other transactions).
When the block size is small relative to the number of transactions,
there is a significant gap between the utilities of the optimal and FIFO-inclusion
blocks. Furthermore, this gap is, unsurprisingly, quite a bit worse than what is predicted by our
bound. For sharp distributions, the gap, shown
in figure~\ref{fig:gaps-heavy}, is significantly larger.
These distributions more closely model situations with competitive MEV opportunities.

\section{Conclusion and future directions}
In this note, we outlined a simple framework to analyze how FIFO ordering
affects welfare in batched systems, such as public blockchains. While this result---that FIFO ordering can
decrease social welfare---is not surprising, our framework formalizes and 
quantifies this intuition. Our framework also suggests a principled way to compare
different inclusion and ordering mechanisms: how closely do they approximate the 
optimal value of the block-building problem? Analyzing this gap for other
mechanisms, along with strengthening the bound~\eqref{eq:main-bound}, both
present avenues for valuable future work.


\bibliographystyle{alpha}
\bibliography{citations.bib}

\newcommand{\etalchar}[1]{$^{#1}$}
\begin{thebibliography}{BCD{\etalchar{+}}19}

\bibitem[BCD{\etalchar{+}}19]{eip1559}
Vitalik Buterin, Eric Conner, Rick Dudley, Matthew Slipper, Ian Norden, and
  Abdelhamid Bakhta.
\newblock Eip-1559: Fee market change for eth 1.0 chain, 2019.

\bibitem[BEKS17]{bezanson2017julia}
Jeff Bezanson, Alan Edelman, Stefan Karpinski, and Viral~B Shah.
\newblock Julia: A fresh approach to numerical computing.
\newblock {\em SIAM review}, 59(1):65--98, 2017.

\bibitem[BPA{\etalchar{+}}21]{distributions.jl}
Mathieu Besançon, Theodore Papamarkou, David Anthoff, Alex Arslan, Simon
  Byrne, Dahua Lin, and John Pearson.
\newblock Distributions.jl: Definition and modeling of probability
  distributions in the juliastats ecosystem.
\newblock {\em Journal of Statistical Software}, 98(16):1--30, 2021.

\bibitem[CFK23]{chitra2023credible}
Tarun Chitra, Matheus~VX Ferreira, and Kshitij Kulkarni.
\newblock Credible, optimal auctions via blockchains.
\newblock {\em arXiv preprint arXiv:2301.12532}, 2023.

\bibitem[CMSZ22]{quickorderfairness}
Christian Cachin, Jovana Mi{\'c}i{\'c}, Nathalie Steinhauer, and Luca Zanolini.
\newblock Quick order fairness.
\newblock In {\em International Conference on Financial Cryptography and Data
  Security}, pages 316--333. Springer, 2022.

\bibitem[CS21]{chung2021foundations}
Hao Chung and Elaine Shi.
\newblock Foundations of transaction fee mechanism design.
\newblock {\em arXiv preprint arXiv:2111.03151}, 2021.

\bibitem[DECA22]{diamandis2022dynamic}
Theo Diamandis, Alex Evans, Tarun Chitra, and Guillermo Angeris.
\newblock Dynamic pricing for non-fungible resources.
\newblock {\em arXiv preprint arXiv:2208.07919}, 2022.

\bibitem[Fla22]{MEV_Suave}
Flashbots.
\newblock The future of mev is suave: Flashbots, Nov 2022.

\bibitem[Fou23]{jito2023}
Jito Foundation.
\newblock Solving the mev problem on solana: A guide for stakers, Feb 2023.

\bibitem[{Gur}23]{gurobi}
{Gurobi Optimization, LLC}.
\newblock {Gurobi Optimizer Reference Manual}, 2023.

\bibitem[KDK22]{aequitas}
Mahimna Kelkar, Soubhik Deb, and Sreeram Kannan.
\newblock Order-fair consensus in the permissionless setting.
\newblock In {\em Proceedings of the 9th ACM on ASIA Public-Key Cryptography
  Workshop}, pages 3--14, 2022.

\bibitem[KDL{\etalchar{+}}21]{themis}
Mahimna Kelkar, Soubhik Deb, Sishan Long, Ari Juels, and Sreeram Kannan.
\newblock Themis: Fast, strong order-fairness in byzantine consensus.
\newblock {\em Cryptology ePrint Archive}, 2021.

\bibitem[KPP04]{kellerer2004knapsack}
Hans Kellerer, Ulrich Pferschy, and David Pisinger.
\newblock {\em Knapsack problems}.
\newblock Springer, 2004.

\bibitem[Kur20]{wendy}
Klaus Kursawe.
\newblock Wendy, the good little fairness widget: Achieving order fairness for
  blockchains.
\newblock In {\em Proceedings of the 2nd ACM Conference on Advances in
  Financial Technologies}, pages 25--36, 2020.

\bibitem[KZGJ20]{kelkar2020order}
Mahimna Kelkar, Fan Zhang, Steven Goldfeder, and Ari Juels.
\newblock Order-fairness for byzantine consensus.
\newblock In {\em Advances in Cryptology--CRYPTO 2020: 40th Annual
  International Cryptology Conference, CRYPTO 2020, Santa Barbara, CA, USA,
  August 17--21, 2020, Proceedings, Part III 40}, pages 451--480. Springer,
  2020.

\bibitem[Rou16]{roughgarden2016twenty}
Tim Roughgarden.
\newblock {\em Twenty lectures on algorithmic game theory}.
\newblock Cambridge University Press, 2016.

\bibitem[sno23]{snoopy2023}
snoopy\_mev.
\newblock \url{https://twitter.com/snoopy_mev/status/1629283898453811200}, Feb
  2023.

\bibitem[VK23]{vafadar2023condorcet}
Mohammad~Amin Vafadar and Majid Khabbazian.
\newblock Condorcet attack against fair transaction ordering.
\newblock {\em arXiv preprint arXiv:2306.15743}, 2023.

\bibitem[ZSC{\etalchar{+}}20]{pompe}
Yunhao Zhang, Srinath Setty, Qi~Chen, Lidong Zhou, and Lorenzo Alvisi.
\newblock Byzantine ordered consensus without byzantine oligarchy.
\newblock In {\em 14th USENIX Symposium on Operating Systems Design and
  Implementation (OSDI 20)}, pages 633--649, 2020.

\end{thebibliography}

\appendix

\section{Distribution definitions} \label{app:dists}

\paragraph{Flat distributions.}
We consider three distributions which yield relatively flat efficiencies as they
have a light (\ie,
sub-exponential or sub-Gaussian) tail. 
These distributions have a finite mean and variance, so
we expect samples to cluster tightly together.
The exponential distribution with parameter $\theta$ has probability density 
function
\[
f(x; \theta) = \frac{1}{\theta} e^{-\frac{x}{\theta}}, \quad x > 0.
\]
The log normal distribution with parameters $\mu$ and $\sigma$ has probability
density function
\[
    f(x; \mu, \sigma) = \frac{1}{x \sqrt{2 \pi \sigma^2}}
\exp \left( - \frac{(\log(x) - \mu)^2}{2 \sigma^2} \right),
\quad x > 0.
\]
The Rayleigh distribution with parameter $\sigma$ has probability density function
\[
    f(x; \sigma) = \frac{x}{\sigma^2} e^{-\frac{x^2}{2 \sigma^2}}, \quad x > 0.
\]

\paragraph{Sharp distributions.}
We consider two utility distributions which have heavy tails.
These distributions have an infinite mean and variance, and we expect the
presence of large outliers when sampling, yielding a relatively sharp
distribution of the efficiencies.
The Levy distribution with parameters $\mu$ and $\sigma$ has probability density
function
\[
    f(x; \mu, \sigma) = \sqrt{\frac{\sigma}{2 \pi (x - \mu)^3}}
    \exp \left( - \frac{\sigma}{2 (x - \mu)} \right), \quad x > \mu.
\]
The Pareto distribution with parameter $\alpha$ has probability density function
\[
    f(x; \alpha) = \frac{\alpha }{x^{\alpha + 1}}, \quad x \ge 1.
\]

\begin{figure}[h!]
    \captionsetup[sub]{font=scriptsize}
    \centering
    \begin{subfigure}[t]{0.45\linewidth}
        \centering
        \includegraphics[width=\columnwidth]{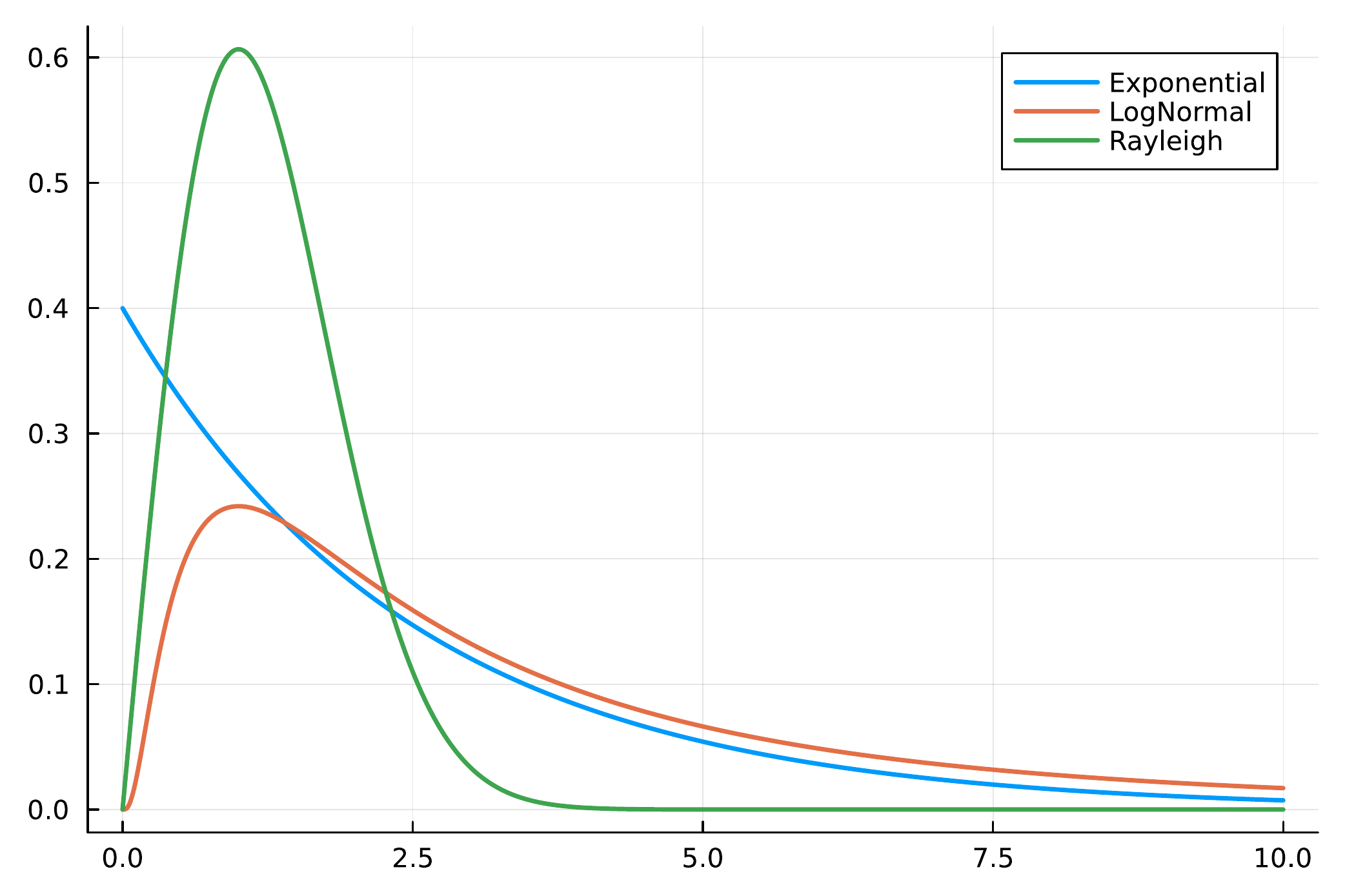}
        \caption{Light-tailed utility distributions.}
    \end{subfigure}
    \hfill
    \begin{subfigure}[t]{0.45\linewidth}
        \centering
        \includegraphics[width=\columnwidth]{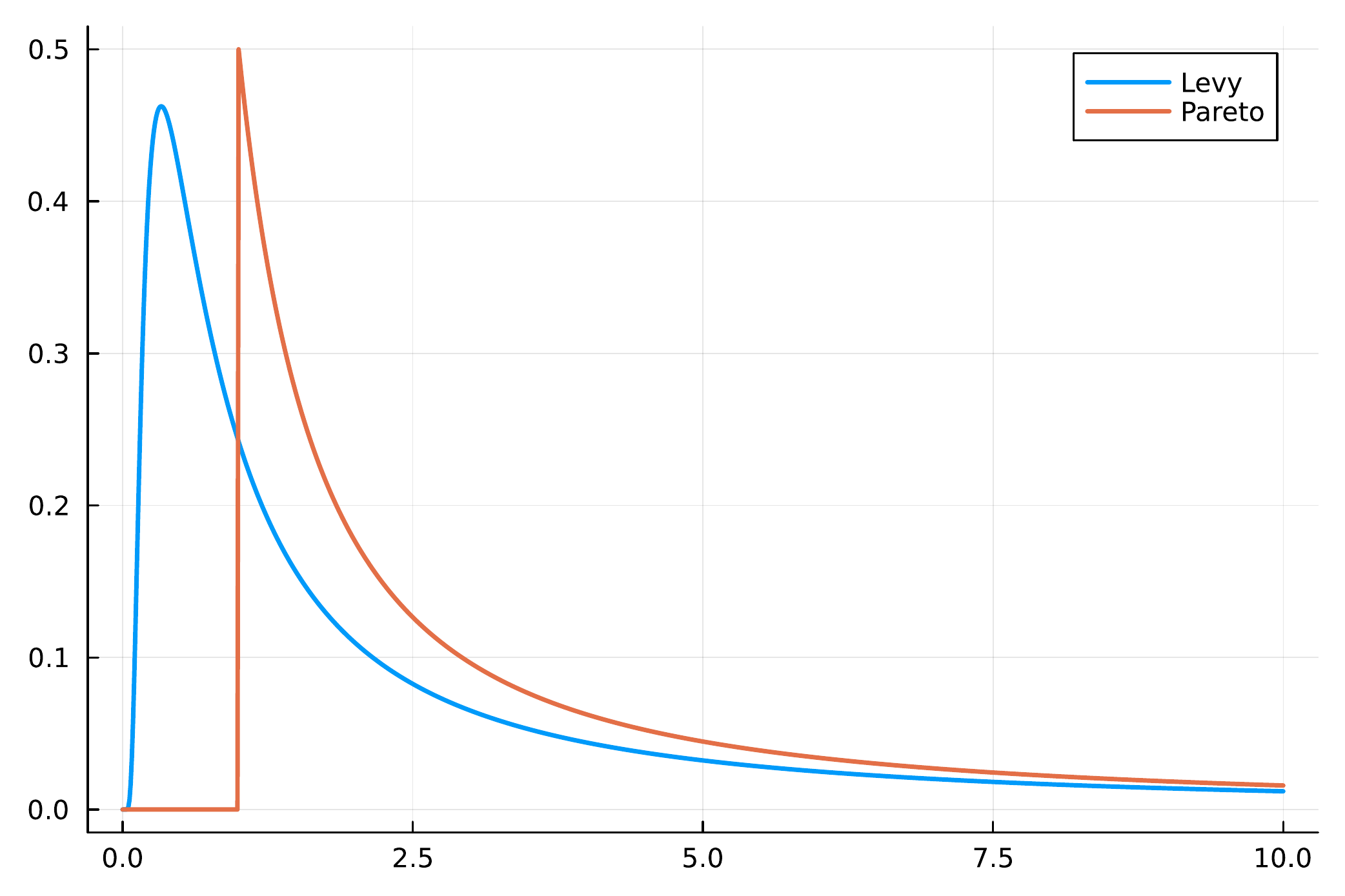}
        \caption{Heavy-tailed utility distributions.}
    \end{subfigure}
    \caption{Distributions used for experiments. Light-tailed utility
        distributions lead to `flat' efficiency distributions, while heavy-tailed
        ones lead to `sharp' efficiency distributions.}
    \label{fig:dists}
\end{figure}

\end{document}